\documentclass[graybox]{svmult}

\usepackage{amsmath}

\usepackage{mathptmx}       
\usepackage{helvet}         
\usepackage{courier}        
\usepackage{type1cm}        
\usepackage{makeidx}         
\usepackage{graphicx}        
\usepackage{multicol}        
\usepackage[bottom]{footmisc}

\makeindex             

\usepackage{a4}
\usepackage{graphicx}
\usepackage{parskip}

\usepackage{amsfonts}

\usepackage{algorithm}
\usepackage{algorithmicx}
\usepackage[noend]{algpseudocode}

\usepackage{natbib}
\usepackage{latexsym}

\newcommand{\R}{\mathbb{R}}

\newcommand{\PP} {{  \rm I\hskip-0.22em P}}

\newcommand{\EE} {{\rm I\hskip-0.48em E}}

\begin{document}

\title*{$\chi^2$-confidence sets in high-dimensional regression}
\author{Sara van de Geer, Benjamin Stucky}
\institute{Sara van de Geer \at Seminar for Statistics, ETH Z\"urich, R\"amistrasse 101,
8092 Z\"urich, Switzerland \email{geer@stat.math.ethz.ch}
}
%
%
\maketitle

\abstract{
We study a high-dimensional regression model. Aim is to construct a confidence
set for a given group of regression coefficients, treating all other regression coefficients as
nuisance parameters. We apply a one-step procedure with the square-root Lasso
as initial estimator and a multivariate square-root Lasso for constructing 
  a surrogate Fisher information matrix.
The multivariate square-root Lasso is based on nuclear norm loss with $\ell_1$-penalty.
We show that this procedure leads to an asymptotically $\chi^2$-distributed
pivot, with a remainder term depending only on the $\ell_1$-error
of the initial estimator. We show
that under $\ell_1$-sparsity conditions on the
regression coefficients $\beta^0$ the square-root Lasso produces to a consistent estimator
of the noise variance and we establish sharp oracle inequalities 
which show that the remainder term is small under further sparsity
conditions on $\beta^0$ and compatibility conditions on the
design.}

\section{Introduction}
Let $X$ be a given $n \times p$ input matrix and $Y$ be a random $n$-vector  of responses.
We consider the high-dimensional situation where the number of variables $p$
exceeds  the number of observations $n$.
 The expectation of $Y$ (assumed to exist) is denoted
by
$f^0 :=\EE Y$. We assume that $X$ has rank $n $ ($n< p $) and
 let $\beta^0$ be any solution of
the equation $X \beta^0 = f^0 $. Our aim is to construct a confidence interval
for a pre-specified group of coefficients $\beta_J^0 := \{ \beta_j^0 : \ j \in J\} $ where $J \subset \{ 1 , \ldots , p \}$ is
a subset of the indices. In other words, the $|J|$-dimensional vector
$\beta_J^0$ is the parameter of interest and all the other coefficients $\beta_{-J}^0 := \{ \beta_j^0 : \ 
j \notin J \} $ are nuisance parameters. 

For one-dimensional parameters of interest ($|J|=1$) the approach in this paper is closely related to earlier work.
The method is introduced in \cite{zhang2014confidence}. Further references are
\cite{jamo13} and 
\cite{van2013asymptotically}. Related approaches can be found in
\cite{Belloni2013b}, \cite{Belloni2013a} and \cite{beletal11}.

For confidence sets for groups of variables ($|J| >1$) one usually would like
to take the dependence between estimators of single parameters into account.
An important paper that carefully does this for confidence sets in $\ell_2$ is
\cite{Mitra}. Our approach is related but differs in an important way.
As in \cite{Mitra} we propose a de-sparsified estimator which is (potentially) asymptotically linear.
However, \cite{Mitra} focus at a remainder term which is small also for large groups.
 Our goal is rather to present a construction which has
a small remainder term after studentizing and which does not rely on strong conditions
on the design $X$. In particular we do not assume any
sparsity conditions on the design. 

The construction involves the square-root Lasso $\hat \beta$
which is introduced by \cite{belloni2011square}. See
Section \ref{square-root.section} for the definition of the estimator $\hat \beta$.
We present a 
a multivariate extension of the square-root Lasso
which takes the nuclear norm of the multivariate residuals as loss function.
Then  we define in Section \ref{construction.section} a de-sparsified estimator $\hat b_J$ of
$\beta_J^0$ which has the form of a one-step estimator with $\hat \beta_J$ as initial estimator and with
multivariate square-root Lasso invoked to obtain a surrogate Fisher information matrix. 
We show that when $Y \sim {\cal N}_n ( f^0 , \sigma_0^2 I)$
(with both $f^0$ and $\sigma_0^2$ unknown), 
a studentized version of $\hat b_J- \beta_J^0$
has asymptotically a $|J|$-dimensional standard normal distribution.

More precisely we will show in Theorem \ref{main.theorem} that for a given 
$|J| \times |J|$ matrix $M=M_{\lambda}$ depending only on $X$ and on a tuning parameter $\lambda$,
one has $ M_{\lambda} (\hat b_J - \beta_J^0 ) /\sigma_0 = {\cal N}_{|J|} (0, I) + {\rm rem} $
where the remainder term ``${\rm rem}$" can be bounded by $
\| {\rm rem} \|_{\infty} \le \sqrt n \lambda \| \hat \beta_{-J} - \beta_{-J}^0 \|_1/\sigma_0$.
The choice of the tuning parameter $\lambda$ is ``free" (and not depending on $\sigma_0$), it can 
can for example be taken of order $\sqrt {\log p / n}$.
The unknown parameter $\sigma_0^2$ can be estimated by the normalized residual sum of
squares $\hat \sigma^2:= \| Y - X \hat \beta \|_2^2/n $ of the square-root Lasso $\hat \beta$.
We show in Lemma \ref{sigma.lemma} that under sparsity conditions on $\beta^0$ one has
$\hat \sigma^2 / \sigma_0^2 = 1+ o_{\PP} (1)$ and then in Theorem \ref{oracle.theorem} an
oracle inequality for the square-root Lasso under 
further sparsity conditions on $\beta^0$ and
compatibility conditions on the design.
The oracle result allows
one to ``verify" when
$\sqrt n \lambda \| \hat \beta_{-J} - \beta_{-J}^0 \|_1/\sigma_0 = o_{\PP} (1)$ so that the
remainder term ${\rm rem}$ is negligible.  An illustration assuming weak sparsity conditions is
given in Section \ref{lr.section}.
As a consequence
$$ \| M_{\lambda} ( \hat b_J - \beta_J^0 ) \|_2^2 / \hat \sigma^2 = \chi_{|J|}^2 ( 1+ o_{\PP }(1)) , $$
where $\chi_{|J|}^2$ is a random variable having a $\chi^2$-distribution with $|J|$ degrees of freedom.
For $|J|$ fixed one can thus construct asymptotic confidence sets for $\beta_J^0$
(we will also consider the case $|J | \rightarrow \infty$ in Section \ref{discussion.section}). 
We however do not control the size of these sets. Larger values for $\lambda$ makes the 
confidence sets
smaller but will also give a larger remainder term. 

In Section \ref{structured.section} we extend the theory to structured
sparsity norms other than $\ell_1$, for example the norm used for the (square-root)
group Lasso, where the demand for $\ell_2$-confidence sets for groups comes
up quite naturally. Section \ref{discussion.section} contains a discussion.
The proofs are in Section \ref{proofs.section}. 

\subsection{Notation}

The mean vector of $Y$ is denoted by $f^0$ and the noise is
$\epsilon := Y - f^0$. 
For a vector $v \in \R^n$ we write (with a slight abuse of notation)
$\| v \|_n^2 := v^T v / n $. We let $\sigma_0^2 := \EE \| \epsilon \|_n^2$
(assumed to exist).

For a vector $\beta \in \R^p$ we set $S_{\beta} := \{ j: \ \beta_j \not= 0 \} $.
For a subset $J \subset \{ 1 , \ldots , p \}$ and
a vector $\beta \in \R^p$ we use the same notation $\beta_J$ for
the $|J|$-dimensional vector $\{ \beta_j : \ j \in J \}$ and the
$p$-dimensional vector $\{ \beta_{j,J} := \beta_j {\rm l} \{ j \in J \} :\ j=1 , \ldots , p \}$. 
The last version allows us to write $\beta = \beta_J + \beta_{-J}$ with
$\beta_{-J} = \beta_{J^c}$, $J^c$ being the complement of the set $J$. 
The $j$-th column of $X$  is denoted by $X_j$ ($j=1 , \ldots , p$). We let $X_J := \{
X_j: \ j \in J \}$ and $X_{-J} := \{ X_j : \ j \notin J \} $. 

For a matrix $A$ we let $\| A \|_{\rm nuclear} := {\rm trace} ( (A^T A)^{1/2} )$
be its nuclear norm. The $\ell_1$-norm of the matrix $A$ is defined as
$\| A \|_1 := \sum_{k} \sum_j | a_{k,j} |$. Its $\ell_{\infty}$-norm is
$\| A \|_{\infty} := \max_k \max_j | a_{k,j} | $.

\section{The square-root Lasso and its multivariate version}\label{square-root.section}

\subsection{The square-root Lasso}
The square-root Lasso (\cite{belloni2011square}) $\hat \beta$ is 
\begin{equation} \label{square-root.equation} \hat \beta := \arg \min_{\beta \in \R^p} \biggl \{ \| Y - X \beta \|_n +
 \lambda_0 \| \beta \|_1 \biggr \} . \end{equation}
 The parameter $\lambda_0>0$ is a tuning parameter. 
 Thus $\hat \beta$ depends on $\lambda_0$ but we do not express this in our notation.
 
 The square-root Lasso can be seen as a method that estimates $\beta^0$ and
the noise variance $\sigma_0^2$ simultaneously. Defining the residuals $\hat \epsilon := Y - X \hat \beta $ 
and letting $\hat \sigma^2 := \| \hat \epsilon \|_n^2 $
one clearly has
\begin{equation} \label{square-root2.equation}
( \hat \beta , \hat \sigma^2 ) = 
\arg \min_{\beta \in \R^p , \ \sigma^2 > 0 } 
\biggl \{ { \| Y- X \beta \|_n^2  \over \sigma } + \sigma + 2 \lambda_0 \| \beta \|_1 \biggr \}  
\end{equation}
provided the minimum is attained at a positive value of $\sigma^2$. 

We note in passing that the square-root Lasso is {\it not} a quasi-likelihood estimator as
the function $ \exp[ -z^2/\sigma - \sigma]$, $z \in \R $, is not a density with respect to
a dominating measure not depending on $\sigma^2>0$. 
The square-root Lasso is moreover not to be confused with the scaled Lasso.
The latter {\it is} a quasi-likelihood estimator. It is studied in
 e.g. \cite{sunzhang11}.

 We show in Section \ref{sigma.section} (Lemmas \ref{Gauss.lemma} and \ref{sigma.lemma})
 that
 for the case where $\epsilon \sim {\cal N}_n (0, \sigma_0^2 I) $ for example
 one has $\hat \sigma \rightarrow \sigma_0$ under $\ell_1$-sparsity conditions on $\beta^0$. In
 Section
 \ref{oracle.section} we 
 establish oracle results for $\hat \beta$ under further sparsity conditions on $\beta^0$ and
 compatibility conditions on $X$ (see Definition \ref{compatibility.definition} for the latter).
 These results hold for a 
``universal" choice of $\lambda_0$ provided an $\ell_1$-sparsity condition on $\beta^0$ is met.

In the proof of our main result in Theorem \ref{main.theorem}, the so-called {\it Karush-Kuhn-Tucker conditions}, or KKT-conditions,
play a major role. Let us briefly discuss these here.
The  KKT-conditions for the square-root Lasso say that
\begin{equation} \label{KKT.equation}
{ X^T (Y- X \hat \beta ) /n \over \hat \sigma} = \lambda_0 \hat z 
\end{equation}
where $\hat z $ is a $p$-dimensional vector with $\| \hat z \|_{\infty} \le 1 $ and
with $\hat z_j = {\rm sign} (\hat \beta_j)$ if $\hat \beta_j \not= 0 $. This
follows 
from sub-differential calculus which defines the sub-differential of
the absolute value function $x \mapsto |x|$ as 
$$\{ {\rm sign} (x) \} \{ x \not= 0 \} + [-1,1] \{ x=0 \} . $$
Indeed, for a fixed $\sigma >0$ the sub-differential with respect to $\beta$ of the expression in curly brackets
given in (\ref{square-root2.equation}) is equal to
$$ -{2 X^T ( Y - X\beta )/ n \over \sigma } + 2 \lambda_0 z (\beta) $$
with, for $j = 1 , \ldots , p$, $z_j (\beta)$ the sub-differential of
$\beta_j \mapsto | \beta_j | $. Setting this to zero at $(\hat \beta , \hat \sigma)$
gives the above KKT-conditions (\ref{KKT.equation}).

\subsection{The multivariate square-root Lasso}\label{multivariate.section}
In our construction of confidence sets we will consider the regression of $X_J$
on $X_{-J}$ invoking a multivariate version
of the square-root Lasso.  To explain the latter, we use here a standard notation
with $X$ being the input and
$Y$ being the response. We will then replace $X$ by $X_{-J}$ and $Y$ by $X_J$
in Section \ref{construction.section}. 

The matrix $X$ is as before an $n \times p$
input matrix and the response $Y$ is now
an $n \times q$ matrix for some $q\ge 1$. We define the multivariate square-root Lasso
\begin{equation} \label{multi-square-root.equation}
 \hat B := \arg \min_{B} \biggr \{\| Y - X B \|_{\rm nuclear}/\sqrt n  + \lambda_0 \| B \|_1 \biggr \} 
 \end{equation} 
 with $\lambda_0 > 0$ again a tuning parameter.  The minimization
 is over all $p \times q$ matrices $B$. 
 We consider $\hat \Sigma := ( Y- X \hat B)^T (Y- X \hat B) / n $ as estimator of the
 noise co-variance matrix. 
 
 The KKT-conditions for the multivariate square-root Lasso will be a major ingredient
 of the proof of the main result in Theorem \ref{main.theorem}. We present
 these KKT-conditions in the following lemma in equation (\ref{multi-KKT.equation}). 
 
 \begin{lemma} \label{KKT.lemma}
We have
 $$ (\hat B , \hat \Sigma ) = \arg \min_{B , \ \Sigma >0 }\biggl \{ 
 {\rm trace} \biggl ((Y- X B)^T (Y-XB) \Sigma^{-1/2} \biggr )/n $$ $$
  \ \ \ \ \ \ \ \ \ \ \ \ \ \ \ \ + {\rm trace} ( \Sigma^{1/2}) + 2 \lambda_0 
 \| B \|_1 \biggr \} $$
 where the minimization is over all symmetric positive definite matrix $\Sigma$
 (this being denoted by $\Sigma > 0$) and  where it is assumed that the minimum is indeed
 attained at some $\Sigma >0$. 
 The multivariate Lasso satisfies the KKT-conditions
  \begin{equation}\label{multi-KKT.equation}
 X^T (Y-X \hat B ) \hat \Sigma^{-1/2} /n = \lambda_0 \hat Z ,
  \end{equation}
  where $\hat Z $ is a $p \times q $ matrix with $\| \hat Z \|_{\infty} \le 1$
  and with $\hat Z_{k,j} = {\rm sign} (\hat B_{k,j} )$ if $\hat B_{k,j} \not= 0 $
  ($k=1 , \ldots , p $, $j=1 , \ldots , q$). 
  \end{lemma}

\section{Confidence sets for $\beta_J^0$}

\subsection{The construction}\label{construction.section}
Let $J \subset \{ 1 , \ldots , p \} $. 
 We are interested in building a confidence set for $\beta_J^0:=
 \{ \beta_j^0 : \ j \in J \} $. 
 To this end, we compute the multivariate ($|J|$-dimensional) square root Lasso
 \begin{equation}\label{gamma.equation}
  \hat \Gamma_J  := \arg \min_{\Gamma_J } \biggl \{ 
 \|  X_J - X_{-J} \Gamma_J  \|_{\rm nuclear} / \sqrt n +
  \lambda  \|\Gamma_J \|_{1}  \biggr \}  
  \end{equation}
  where $\lambda >0$ is a tuning parameter.
 The minimization is over all $(p-|J|) \times |J|$ matrices $\Gamma_J$.
 We let
  \begin{equation}\label{tildeT.equation}
 \tilde T_J :=  ( X_J - X_{-J} \hat \Gamma_J )^T X_J/ n 
 \end{equation}
 and
 \begin{equation}\label{hatT.equation}
 \hat T_J := ( X_J - X_{-J} \hat \Gamma_J )^T  ( X_J - X_{-J} \hat \Gamma_J )/ n ,
 \end{equation}
 We assume throughout that the ``hat" matrix
 $\hat T_J$ is non-singular. The ``tilde" matrix $\tilde T_J$ only needs to be non-singular
 in order that the de-sparsified estimator $\hat b_J$ given below in Definition \ref{de-sparse.definition}
 is well-defined. However, for the normalized version we need not assume non-singularity of $\tilde T_J$.

 The KKT-conditions (\ref{multi-KKT.equation}) appear in the form
 \begin{equation}\label{gammaKKT.equation}
 X_{-J}^T ( X_J - X_{-J} \hat \Gamma_J ) \hat T_J^{-1/2}/n = \lambda \hat Z_J , 
 \end{equation}
 where $\hat Z_J$ is a $(p-|J|) \times |J|$ matrix with $(\hat Z_J)_{k,j} = {\rm sign} (\hat \Gamma_J) _{k,j}$
 if $(\hat \Gamma_J) _{k,j}\not= 0 $ and $\| \hat Z_J \|_{\infty} \le 1 $.

 We define the normalization matrix 
 \begin{equation}\label{M.equation}
 M := M_{\lambda} := \sqrt n \hat T_J^{-1/2} \tilde T_J . 
 \end{equation}
 
 \begin{definition}\label{de-sparse.definition}
 The {\rm de-sparsified} estimator of $\beta_J^0$ is
 $$ \hat b_J := \hat \beta_J + \tilde T_J^{-1} (X_J - X_{-J} \hat \Gamma_J)^T (Y - X \hat \beta) /n , $$
 with $\hat \beta$ the square-root Lasso given in (\ref{square-root.equation}),
 $\hat \Gamma_J$ the multivariate square-root Lasso given in (\ref{gamma.equation}) and
 the matrix $\tilde T_J$ given in (\ref{tildeT.equation}).
 The normalized de-sparsified estimator is
 $M \hat b_J   $
 with $M$ the normalization matrix given in (\ref{M.equation}). 
 \end{definition}

 \subsection{The main result} \label{main.section}
 
 Our main result is rather simple. It shows that using the
 multivariate square-root Lasso for de-sparsifying, and then
 normalizing, results in a well-scaled ``asymptotic pivot" (up to the estimation of
 $\sigma_0$ which we will do in the next section). Theorem
 \ref{main.theorem} actually does not require $\hat \beta$ to be the square-root
 Lasso but for definiteness we have made this specific choice throughout the paper
 (except for Section \ref{structured.section}).

 \begin{theorem} \label{main.theorem} Consider the model $Y \sim {\cal N}_n ( f^0 , \sigma_0^2)$
 where $f^0 = X \beta^0$. 
 Let $\hat b_J$ be the de-sparsified estimator given
 in Definition \ref{de-sparse.definition} and let $M \hat b_J$ be its
 normalized version. Then
 $$M( \hat b_J - \beta_J^0 ) /\sigma_0 = {\cal N}_{|J|} (0 , I) + {\rm rem} $$
 where $\| {\rm rem} \|_{\infty} \le \sqrt n \lambda \| \hat \beta_{-J} - \beta_{-J}^0 \|_1 / \sigma_0 $. 
 
 \end{theorem}

To make Theorem \ref{main.theorem} work we need to bound 
$\| \hat \beta - \beta^0 \|_1 / \hat \sigma$ where $\hat \sigma$ is an estimator of $\sigma_0$.
This is done in Theorem \ref{oracle.theorem} with $\hat \sigma$ the estimator $\| \hat \epsilon \|_n$
from the square-root Lasso.
A special case is presented in Lemma \ref{bound.lemma} which imposes
weak sparsity conditions for $\beta^0$. Bounds for
$\sigma_0 / \hat \sigma$ are also given.

Theorem \ref{main.theorem} is about the case where the noise $\epsilon$ is i.i.d.\ normally distributed.
 This can be generalized as from the proof 
 we see that the ``main" term is linear in $\epsilon$. For independent errors
 with common variance $\sigma_0^2$ say, one needs to assume
 the Lindeberg condition for establishing asymptotic normality.

\section{Theory for the square root Lasso}\label{theory.section}
Let $f^0 := \EE Y$, $\beta^0$ be a solution of $X \beta^0 = f^0$ and
define $\epsilon := Y- f^0$.
Recall the square-root Lasso $\hat \beta$ given in (\ref{square-root.equation}).
It depends on the tuning parameter $\lambda_0>0$.
In this section we develop theoretical bounds, which are closely
related to results in \cite{sun2013sparse} (who by the way use
the term scaled Lasso instead of square-root Lasso in that paper).
There are two differences. Firstly, our lower bound for the residual sum of squares
of the square-root Lasso requires, for the case where no conditions are imposed on the compatibility constants,
a smaller value for the tuning parameter (see Lemma \ref{sigma.lemma}). 
 These compatibility constants, given in Definition \ref{compatibility.definition}, are required only later for
 the oracle results.
 Secondly,
we establish an oracle inequality that is sharp (see Theorem \ref{oracle.theorem} in Section
\ref{oracle.section} where we present more details).

Write $\hat \epsilon:= Y - X \hat \beta$ 
 and $\hat \sigma^2 := \| \hat \epsilon \|_n^2$. We consider bounds
in terms of $\| \epsilon \|_n^2$, the ``empirical" variance of the unobservable noise.
This is a random quantity but under obvious conditions it convergences
to its expectation $\sigma_0^2$. 
Another random quantity that appears in our bounds is $\epsilon / \| \epsilon \|_n$, which is
a random point on the $n$-dimensional unit sphere.
We write
$$ \hat R:= { \| X^T \epsilon \|_{\infty}  \over n \| \epsilon \|_n } .$$
When all $X_j$ are normalized such that $\| X_j \|_n=1$, the quantity $\hat R$
is the maximal ``empirical" correlation between noise and input variables.
Under distributional assumptions $\hat R$ can be bounded with large probability by some
constant $R$.
For completeness we work out the case of i.i.d.\ normally distributed errors.

\begin{lemma} \label{Gauss.lemma} Let $\epsilon \sim {\cal N}_n (0 , \sigma_0^2 I)$. 
Suppose the normalized case where $\| X_j \|_n =1$ for all $j=1 , \ldots , p $. 
Let $\alpha_0$, $\underline \alpha$ and $\bar \alpha$ be given positive error levels
such that $\alpha_0 + \underline \alpha + \bar \alpha < 1$ and
$\log(1/ {\underline \alpha} ) < n/4$. 
Define 
$$ {\underline \sigma }^2: = \sigma_0^2 \biggl ( 1- 2 \sqrt {\log (1/ {\underline \alpha}) \over n} \biggr ) ,$$
$$ {\bar \sigma}^2 := \sigma_0^2 \biggl ( 1+ 2 \sqrt {\log (1/ {\bar \alpha}) \over n} +   {2 \log (1/ {\bar \alpha}) \over n} 
\biggr ) $$
and
$$R:= \sqrt {   \log (2p/\alpha_0)   \over n- 2 \sqrt {n \log (1/ {\underline \alpha})  } } . $$
We have
$$\PP ( \| \epsilon \|_n \le {\underline \sigma}  ) \le {\underline \alpha}, \  \PP ( \| \epsilon \|_n \ge \bar \sigma) \le \bar \alpha  $$
and
$$\PP  ( \hat R \ge  R \ \cup \| \epsilon \|_n \le {\underline \sigma} ) \le \alpha_0 + {\underline \alpha}  . $$
\end{lemma}

\subsection{Preliminary lower and upper bounds for $\hat \sigma^2$}
\label{sigma.section}

   We now show that the estimator of the variance $\hat \sigma^2= \| \hat \epsilon \|_n^2$, obtained by applying the
   square-root Lasso, converges to the noise variance $\sigma_0^2$. The result holds without conditions
   on compatibility constants (given in Definition \ref{compatibility.definition}). We do however need
   the $\ell_1$-sparsity condition (\ref{eta.equation}) on $\beta_0$. This condition will be discussed below
   in an asymptotic setup.

    \begin{lemma}\label{sigma.lemma}
    Suppose that for some $0 < \eta < 1 $, some $R>0$ and some ${\underline \sigma } > 0$, we have
    $$\lambda_0 (1- \eta)  \ge R $$
    and
    \begin{equation}\label{eta.equation}
     \lambda_0 \| \beta^0 \|_1 / {\underline \sigma}  \le 2 \biggl ( \sqrt { 1 + (\eta/2)^2} -1 \biggr ) . 
     \end{equation}
    Then on the set where $\hat R \le R$ and $\| \epsilon \|_n \ge {\underline \sigma}$ we have
  $\biggl | \| \hat \epsilon \|_n/ \| \epsilon \|_n - 1 \biggr | \le \eta  $.
    
    \end{lemma} 
    
    We remark here that the  the
   result of Lemma \ref{sigma.lemma} is also useful when using a square-root Lasso for constructing an asymptotic confidence
   interval for a single parameter, say $\beta_j^0$. Assuming random design it can be applied to
    show that without imposing compatibility conditions the
   residual variance of the square root Lasso for the regression of $X_j$
   on all other variables $X_{-j}$ does not degenerate.

    {\bf Asymptotics}
     Suppose $\epsilon_1 , \ldots , \epsilon_n$ are i.i.d. with finite variance $\sigma_0^2$. 
     Then clearly $\| \epsilon \|_n /\sigma_0 \rightarrow 1$ in probability. 
     The normalization in (\ref{eta.equation}) by $\underline \sigma$ -
     which can be taken more or less equal to $\sigma_0$ - makes sense if we
     think of the standardized model
     $$ \tilde Y = X \tilde \beta^0 + \tilde \epsilon, $$
     with $\tilde Y = Y / \sigma_0$, $\tilde \beta^0 = \beta^0 / \sigma_0$ and
     $ \tilde \epsilon = \epsilon / \sigma_0$. The condition (\ref{eta.equation}) is
     a condition on the normalized $\tilde \beta^0$.
     The rate of growth assumed there is quite common.
     First of all, it is clear that if $\| \beta^0 \|_1$ is very large then the estimator is
     not  very good because of the penalty on large values of $\| \cdot \|_1$.
     The condition (\ref{eta.equation}) is moreover closely related to standard
     assumptions in compressed sensing.
     To explain this we first note that 
    $$\| \beta^0 \|_1 / \sigma_0 \le   \sqrt{ s_0}  \|  \beta^0 \|_2 /\sigma_0$$
    when $s_0$ is the number of non-zero entries of $\beta^0$
     (observe that $s_0$ is a scale free property of $\beta^0$).
     The term $\| \beta^0 \|_2 / \sigma_0$ can be seen as a signal-to-noise ratio.
     Let us assume this signal-to-noise
    ratio stays bounded. 
     If $\lambda_0$ corresponds to the standard choice $\lambda_0 \asymp \sqrt {\log p / n }$ 
    the assumption (\ref{eta.equation}) holds with $\eta= o (1)$ as soon as we assume 
     the standard assumption $s_0 =o( n/  \log p )$.

    \subsection{An oracle inequality for the square-root Lasso}\label{oracle.section}
    
    Our next result is an oracle inequality for the square-root Lasso. It is
    as the corresponding result for the Lasso as established in  \cite{bickel2009sal}.
    The oracle inequality of Theorem \ref{oracle.theorem} is sharp in the
    sense that there is a constant 1 in front of the approximation error
    $\| X (\beta - \beta^0) \|_n^2$ in (\ref{oracle.equation}). This sharpness is obtained along the lines of arguments from
    \cite{koltchinskii2011nuclear}, who prove sharp oracle inequalities for the Lasso and for
    matrix problems. We further
    have extended the situation in order to establish an oracle inequality for
    the $\ell_1$-estimation error $\| \hat \beta - \beta^0 \|_1$ where we use arguments
    from \cite{geer2014weakly} for the Lasso.
    For the square-root Lasso, the paper \cite{sun2013sparse} also has oracle inequalities, but
    these are not sharp. 
    
    Compatibility constants are introduced in \cite{vandeGeer:07a}. They play a role in the
    identifiability of $\beta^0$. 
    
    \begin{definition} \label{compatibility.definition} Let $L>0$ and $S \subset \{ 1 , \ldots , p \}$.
    The compatibility constant is
    $$ \hat \phi^2 ( L, S) = \min \biggl \{ |S| \| X \beta \|_n^2 :\ \| \beta_S \|_1 = 1 , \ \| \beta_{-S} \|_1 \le L \biggr \} . $$
   \end{definition}
   
   We recall the notation $S_{\beta} = \{ j: \ \beta_j \not= 0 \}$
    appearing in (\ref{oracle.equation}).
    
          \begin{theorem}\label{oracle.theorem} 
Let $\lambda_0$ satisfy for some $R>0$
$$\lambda_0 (1- \eta) > R  $$ 
and assume the $\ell_1$-sparsity (\ref{eta.equation}) for some $0 < \eta < 1$ and $\underline \sigma >0$, i.e.
$$\lambda_0 \| \beta^0 \|_1 / {\underline \sigma} \le 2 \biggl ( \sqrt { 1 + (\eta/2)^2} -1 \biggr ) . $$
Let $0 \le \delta < 1$ be arbitrary and define
$$ \underline {\lambda} := \lambda_0(1-\eta) - R ,$$ 
$$\bar \lambda := \lambda_0(1+ \eta) +R  + \delta \underline{\lambda} $$
and
$$L:=  { \bar \lambda \over (1- \delta ) \underline{\lambda}  } . $$
Then on the set where $\hat R \le R$ and $\| \epsilon \|_n \ge {\underline \sigma} $, we have
$$ 2\delta  \underline{\lambda}  \| \hat \beta - \beta^0 \|_1\| \epsilon \|_n +  \| X ( \hat \beta - \beta^0 ) \|_n^2 
 $$ 
 \begin{equation} \label{oracle.equation}
\le  \min_{S} \biggl \{ \min_{\beta \in \R^p, \ S_{\beta} = S } \biggl [2\delta \underline{\lambda} \| \beta - \beta^0 \|_1
 \| \epsilon \|_n  + 
 \| X( \beta - \beta^0) \|_n^2    \biggr ]  
+  {\bar \lambda}^2 { |S | \| \epsilon \|_n^2 \over  \hat \phi^2 (L,S) }\biggr \} .  
\end{equation}
\end{theorem}

  The result of Theorem \ref{oracle.theorem} leads to a trade-off between the approximation error
  $\| X ( \beta - \beta^0) \|_n^2$, the $\ell_1$-error $\| \beta - \beta^0 \|_1$ and the sparseness\footnote{or non-sparseness actually} 
  $|S_{\beta}|$ (or rather the {\it effective sparseness}  $|S_{\beta}|/ \hat \phi^2 (L,S_{\beta} ) $).
    
    \section{A bound for the $\ell_1$-estimation error under (weak) sparsity}\label{lr.section}
    In this section we assume 
    \begin{equation}\label{lr.equation}
    \sum_{j=1}^r | \beta_j^0 |^r \le \rho_r^r , 
    \end{equation}
    where $0<  r < 1$ and where $\rho_r > 0$ is a constant that
    is ``not too large".  This
    is sometimes called {\it weak sparsity} as opposed to
    {\it strong sparsity} which requires
    "not too many" non-zero coefficients
       $ s_0 := 
     \# \{ \beta_j^0  \not= 0 \} $.               
    We start with bounding the right hand side of the oracle inequality 
    (\ref{oracle.equation}) in Theorem \ref{oracle.theorem}.
    
     We let $S_0:= S_{\beta^0} $ be the active set $S_0 := \{ j : \beta_j^0 \not= 0 \} $
  of $\beta_0$ and  let $\hat \Lambda_{\rm max}^2 (S_0)$ be the
  largest eigenvalue of $X_{S_0}^T X_{S_0} / n $. 
  The cardinality of $S_0$ is denoted by $s_0 = |S_0 |$.
  We assume in this section the normalization  $\| X_j \|_n =1$ so that
  that $\hat \Lambda_{\rm max} (S_0) \ge 1$ and
   ans $\hat \phi(L, S) \le 1$ for any $L$ and $S$.
    
    \begin{lemma} \label{lr.lemma} Suppose $\beta^0$ satisfies the weak sparsity condition
    (\ref{lr.equation}) for some $0 < r < 1$ and $\rho_r >0$.
   For any positive $\delta $, ${\underline \lambda}$,
  $\bar \lambda $ and $L$
  $$
  \min_{S} \biggl \{ \min_{\beta \in \R^p, \ S_{\beta} = S } \biggl [2\delta \underline{\lambda} \| \beta - \beta^0 \|_1
 \| \epsilon \|_n  + 
 \| X( \beta - \beta^0) \|_n^2    \biggr ]  
+  {\bar \lambda}^2 { |S | \| \epsilon \|_n^2 \over  \hat \phi^2 (L,S) }\biggr \}  $$
$$ \le 
 2 {\bar \lambda}^{2-r} \biggl (  \delta {\underline \lambda}/\bar \lambda    +
 { \hat \Lambda_{\rm max}^r (S_0) \over  \hat \phi^{2}  (L, \hat S_*)}
\biggr ) \biggr ( { \rho_r  \over \| \epsilon \|_n  } \biggr )^r \| \epsilon \|_n^2 ,$$
where ${ \hat S_*} := \{ j : \ | \beta_j^0 | > \bar \lambda \| \epsilon \|_n / \hat \Lambda_{\rm max}  (S_0) \} $. 

\end{lemma}

As a consequence, we obtain  bounds 
for the prediction error and  $\ell_1$-error of the square-root Lasso under (weak) sparsity.
We only present the bound for the $\ell_1$-error as this is what we need
in Theorem \ref{main.theorem} for the construction of asymptotic confidence sets.

To avoid being taken away by all the constants, we make some arbitrary choices
in Lemma \ref{bound.lemma}:
we set $\eta \le 1/3$ in the $\ell_1$-sparsity condition (\ref{eta.equation}) and we set
$\lambda_0 (1- \eta) = 2R$. We choose $\delta = 1/7$.

We include the confidence statements that are given
in Lemma \ref{Gauss.lemma} to complete the picture. 

 \begin{lemma}\label{bound.lemma} 
 Suppose $\epsilon \sim {\cal N}_n (0, \sigma_0^2 I )$. 
 Let $\alpha_0$ and $\underline \alpha$ be given positive error levels
such that $\alpha_0 + \underline \alpha  < 1$ and
$\log(1/ {\underline \alpha} ) < n/4$. 
Define 
$$ {\underline \sigma }^2: = \sigma_0^2 \biggl ( 1- 2 \sqrt {\log (1/ {\underline \alpha}) \over n} \biggr ) , \ 
R:= \sqrt {   \log (2p/\alpha_0)   \over n- 2 \sqrt {n \log (1/ {\underline \alpha})  } } . $$
 Assume the $\ell_1$-sparsity condition
 $$ R\| \beta^0 \|_1 / {\underline \sigma} \le (1-\eta)  \biggl ( \sqrt { 1 + (\eta/2)^2} -1 \biggr ) , \ where \ 0<\eta \le 1/3 $$
 and the $\ell_r$-sparsity condition
    (\ref{lr.equation}) for some $0 < r < 1$ and $\rho_r >0$. 
    Set 
    $${ S_*} := \{ j : \ | \beta_j^0 | >3 R \underline \sigma  / \hat \Lambda_{\rm max} (S_0)   \} . $$
Then for $\lambda_0 (1- \eta)=2R$, 
with probability at least $1- \alpha_0 - {\underline \alpha}$ we have the $\ell_r$-sparsity based bound
$$
    (1- \eta) { \| \hat \beta - \beta^0 \|_1 \over \hat \sigma  }  \le  { \| \hat \beta - \beta^0 \|_1 \over \| \epsilon \|_n } \le 
     {(6 R)}^{1-r} \biggl (  1   +
 { 6^2 \hat \Lambda_{\rm max}^r (S_0) \over  \hat \phi^{2}  (6, S_*)}
\biggr ) \biggr ( { \rho_r  \over {\underline \sigma}  } \biggr )^r ,$$
the $\ell_0$-sparsity based bound
$$ (1- \eta)  { \| \hat \beta - \beta^0 \|_1 \over \hat \sigma  }  \le  { \| \hat \beta - \beta^0 \|_1 \over \| \epsilon \|_n } \le
3R \biggl ( { 6^2 s_0 \over \hat \phi^2 ( 6, S_0) } \biggr ) $$
and moreover the following lower bound for the estimator $\hat \sigma$ of the noise level:
$$ (1- \eta) \sigma_0 / \hat \sigma \le   \biggl (1- 2 \sqrt { \log (1/{\underline \alpha} ) \over n } \biggr )^{-1/2} .$$

 \end{lemma}

{\bf Asymptotics} Application of Theorem \ref{main.theorem} 
with $\sigma_0$ estimated by $\hat \sigma$ requires that $\sqrt n \lambda \| \hat \beta - \beta^0 \|_1 / 
\hat \sigma $  tends to zero in probability.  Taking $\lambda \asymp \sqrt {\log p / n} $ and 
for example $\alpha_0 = \alpha = 1/p$, we see that this is the case under the conditions of
Lemma \ref{bound.lemma} as soon as for some $0< r < 1$ the following $\ell_r$-sparsity based bound holds:
$$ \biggl ( {\hat \Lambda_{\rm max}^r (S_0) \over  \hat \phi^{2}  (6, S_*) } \biggr ) \biggl ({\rho_r \over \sigma_0 } \biggr )^r=   { o(n/ \log p )^{1- r \over 2 } \over {(\log p)^{1 \over 2} }}  .$$
Alternatively,  one may require the $\ell_0$-sparsity based bound
$$ \biggl ( {1 \over  \hat \phi^{2}  (6, S_0) } \biggr ) {s_0  }  = {o(n/ \log p )^{1 \over 2 } \over(\log p)^{1 \over 2} } . $$

%
%
\section{Structured sparsity}\label{structured.section} 
We will now show that the results hold for norm-penalized estimators with norms
other than $\ell_1$.
Let
$\Omega$ be some norm on $\R^{p- |J|}$ and define for a 
$(p- | J|) \times |J|$ matrix $A:= (a_1 , \ldots , a_{|J|})$ 
$$ \| A \|_{1 , \Omega} := \sum_{j=1}^{|J| } \Omega ( a_j ) . $$
For a vector $z \in \R^{p - |J|}$ we define the dual norm
$$\Omega_* (z) = \sup_{\Omega (a) \le 1 } |z^T a | , $$
and for a $(p- | J|) \times |J|$ matrix $Z= (z_1 , \ldots , z_{|J|})$ we let
$$\| Z \|_{\infty, \Omega_*} = \max_{1 \le j \le |J| } \Omega_* ( z_j ) . $$
Thus, when $\Omega$ is the $\ell_1$-norm we have $\| A \|_{1 , \Omega} = \| A \|_1 $ and
$\| Z \|_{\infty , \Omega_*} = \| Z \|_{\infty} $.
We let the multivariate 
square-root $\Omega$-sparse estimator be
$$ \hat \Gamma_J :=
\arg \min_{\Gamma_J} \biggl \{ \| X_J - X_{-J} \Gamma_J \|_{\rm nuclear} / \sqrt n + \lambda
\| \Gamma \|_{1, \Omega} \biggr \} . $$
This estimator equals (\ref{gamma.equation}) when $\Omega$ is the $\ell_1$-norm.

 We let, as in (\ref{tildeT.equation}), (\ref{hatT.equation}) and (\ref{M.equation}) but now
 with the new $\hat \Gamma_J$, the quantities $\tilde T_J$, $\hat T_J$ and $M$ be defined as
  $$
 \tilde T_J :=  ( X_J - X_{-J} \hat \Gamma_J )^T X_J/ n ,
$$
 $$
 \hat T_J := ( X_J - X_{-J} \hat \Gamma_J )^T  ( X_J - X_{-J} \hat \Gamma_J )/ n 
$$
and
$$
 M := M_{\lambda} := \sqrt n \hat T_J^{-1/2} \tilde T_J . 
 $$
 
The $\Omega$-de-sparsified  estimator of $\beta_J^0$ is as in Definition \ref
{de-sparse.definition}
 $$ \hat b_J := \hat \beta_J + \tilde T_J^{-1} (X_J - X_{-J} \hat \Gamma_J)^T (Y - X \hat \beta) /n , $$
 but now with $\hat \beta$ not necessarily the square root Lasso but
 a suitably chosen initial estimator and with
 $\hat \Gamma_J$ the multivariate square-root $\Omega$-sparse estimator.
 The normalized de-sparsified estimator is
 $M \hat b_J   $
 with normalization matrix $M$ given above. We can then easily derive the following
 extension of Theorem \ref{main.theorem}. 
 
\begin{theorem}\label{structured.theorem}
Consider the model $Y \sim {\cal N}_n ( f^0 , \sigma_0^2)$
 where $f^0 = X \beta^0$. 
 Let $\hat b_J$ be the $\Omega$-de-sparsified estimator  depending on some initial estimator
 $\hat \beta$. 
Let $M \hat b_J$ be its
 normalized version. Then
 $$M( \hat b_J - \beta_J^0 ) /\sigma_0 = {\cal N}_{|J|} (0 , I) + {\rm rem} $$
 where $\| {\rm rem} \|_{\infty} \le \sqrt n \lambda \Omega( \hat \beta_{-J} - \beta_{-J}^0 )  / \sigma_0 $. 
 \end{theorem}

 We see from Theorem \ref{structured.theorem} that confidence sets follow
 from fast rates of convergence of the $\Omega$-estimation error.
 The latter is studied in \cite{Bach:10},
\cite{Obozinski:12}
and \cite{geer2014weakly} for the case where  the initial estimator is the least
squares estimator with penalty based on a sparsity inducing norm $\bar \Omega$ (say). Group sparsity
\cite{lin06grouplasso} is an example which we shall now briefly discuss.

\begin{example} \label{group.example}Let $G_1 , \ldots , G_T$ be given mutually
disjoint subsets of $\{ 1 , \ldots , p \}$ and take as sparsity-inducing norm
$$ \bar \Omega (\beta) := \sum_{t=1}^T \sqrt { | G_t | } \| X\beta_{G_t} \|_2, \ \beta \in \R^p . $$
The group Lasso is the minimizer of least squares loss with penalty proportional to $\bar \Omega$.
Oracle inequalities for the $\bar \Omega$-error of the group Lasso have been derived in
\cite{lounici:11} for example. For the square-root version we refer to \cite{Bunea13}. With group
sparsity, it lies at hand to consider confidence sets for one of the groups $G_t$ i.e., to take
$J=G_{t_0}$ for a given $t_0$. Choosing
$$\Omega ( a) = \sum_{t \not=t_0} \sqrt { | G_t | } \| Xa_{G_t} \|_2, \ a \in \R^{p-|G_{t_0} | } $$ 
will ensure that $
\Omega ( \hat \beta_{-G_{t_0}} - \beta_{-G_{t_0}} )\le  \bar \Omega (\hat \beta - \beta^0) $
which gives one a handle to control the remainder term in Theorem \ref{structured.theorem}. 
This choice of $\Omega$ for constructing the confidence set makes sense if one
believes that the group structure describing the relation between the response $Y$ and
the input $X$ is also present in the relation between $X_{G_{t_0}}$ and $X_{-G_{t_0}}$. 

\end{example}

\section{Simulations}
Here we denote by $\widetilde{Y}=\widetilde{X}B_{0}+\epsilon$ an arbitrary linear multivariate regression.
In a similar fashion to the suqare-root Algorithm in Buena et al. (2014) we propose the following algorithm for the multivariate square-root Lasso:


\begin{algorithm}
  \caption{msrL}
  \begin{algorithmic}[1]
    \Require{Take a constant $K$ big enough, and choose an arbitrary starting matrix $B(0)\in \mathbb{R}^{p\times q}$.}
    \State $\widetilde{Y} \gets \widetilde{Y}/K$
    \State $\widetilde{X} \gets \widetilde{X}/K$
      \For{$t = 0,1,2,... t_{stop}$}
        \State $B(t+1):=\overline{\Phi}\left(B(t)+X^{T}\cdot (Y-XB(t)); \text{ }\lambda \|Y-XB(t)\|_{nuclear} \right)$
      \EndFor
      \State \Return{$B(t_{stop}+1)$}
  \end{algorithmic}
\end{algorithm}

Here we denote
$$\overline{\Phi}(a;\eta) := \begin{cases} 0, & \mbox{if } a=0 \\ \frac{a}{\|a\|_{2}}(\|a\|_{2}-\eta)_{+}, & \mbox{if } a>0 \end{cases}.$$
The value $t_{stop}$ can be chosen in such a way that one gets the desired accuracy for the algorithm. This algorithm is based on a Fixpoint equation from the KKT conditions.
The square root Lasso is calculated via the algorithm in Buena et al. (2014).

We consider the usual linear regression model:
$$Y=X\beta+\epsilon.$$
In our simulations we take a design matrix X, where the rows are fixed i.i.d. realizations from $\mathcal{N}(0,\Sigma)$. We have $n$ observations, and $p$ explanatory variables.
The covariance matrix $\Sigma$ has the following toeplitz structure $\Sigma_{i,j}=0.9^{|i-j|}$.
The errors are i.i.d. Gaussian distributied, with variance $\sigma^{2}=1$. 
A set of points $J$ is also chosen. $J$ denotes the set of indices of the parameter vector $\beta$ that we want to find asymptotic group confidence intervals for. Define as $q=|J|$ the number of
indices of interest.
For each different setting of $p$ and $n$ we do $r=1000$ simulations.
In each repetition we calculate the teststatistic $\chi^{2}$. A significance level of $0.05$ is chosen.
The lambda of the square root LASSO $\lambda_{srL}$ in the simulations is the theoretical lambda $\lambda_{srLt}$ scaled by $3$. 
For the lambda of the multivariate square root LASSO
$\lambda_{msrL}$ we do not have theoretical results yet. That is why we use cross-validation where we minimize the error expressed in nuclear norm, to define $\lambda_{msrL}$. 
It is important to note, that the choice of $\lambda_{msrL}$ is very crucial, especially
for cases where $n$ is small. One could tune $\lambda_{msrL}$ in such a fashion that even cases like $n=100$ work much better, see figure \ref{fig2}. 
But the point here is to see what happens to the chi-squared test statistic with a fixed rule for the choice of $\lambda_{msrL}$.
This basic set up is used throughout all the simulations below.

\subsection{Asymptotic distribution}
\begin{figure}[h!t]

  \centering
    \includegraphics[width=1\textwidth]{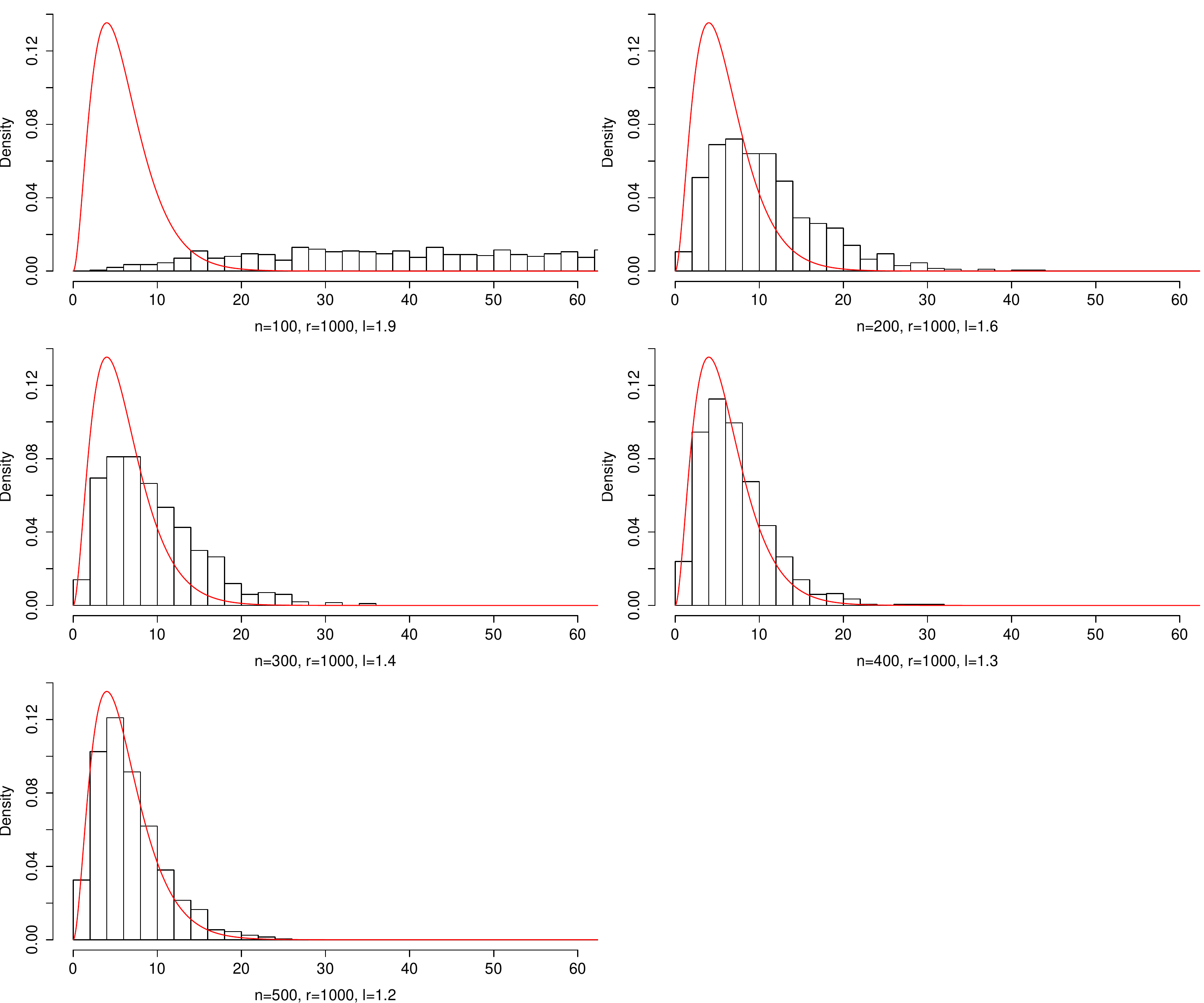}
      \caption{Histogram of Teststatistic, l=$\lambda$, the cross-validation lambda}
      \label{fig1}
\end{figure}

\begin{figure}[h!t]

    \includegraphics[width=0.5\textwidth]{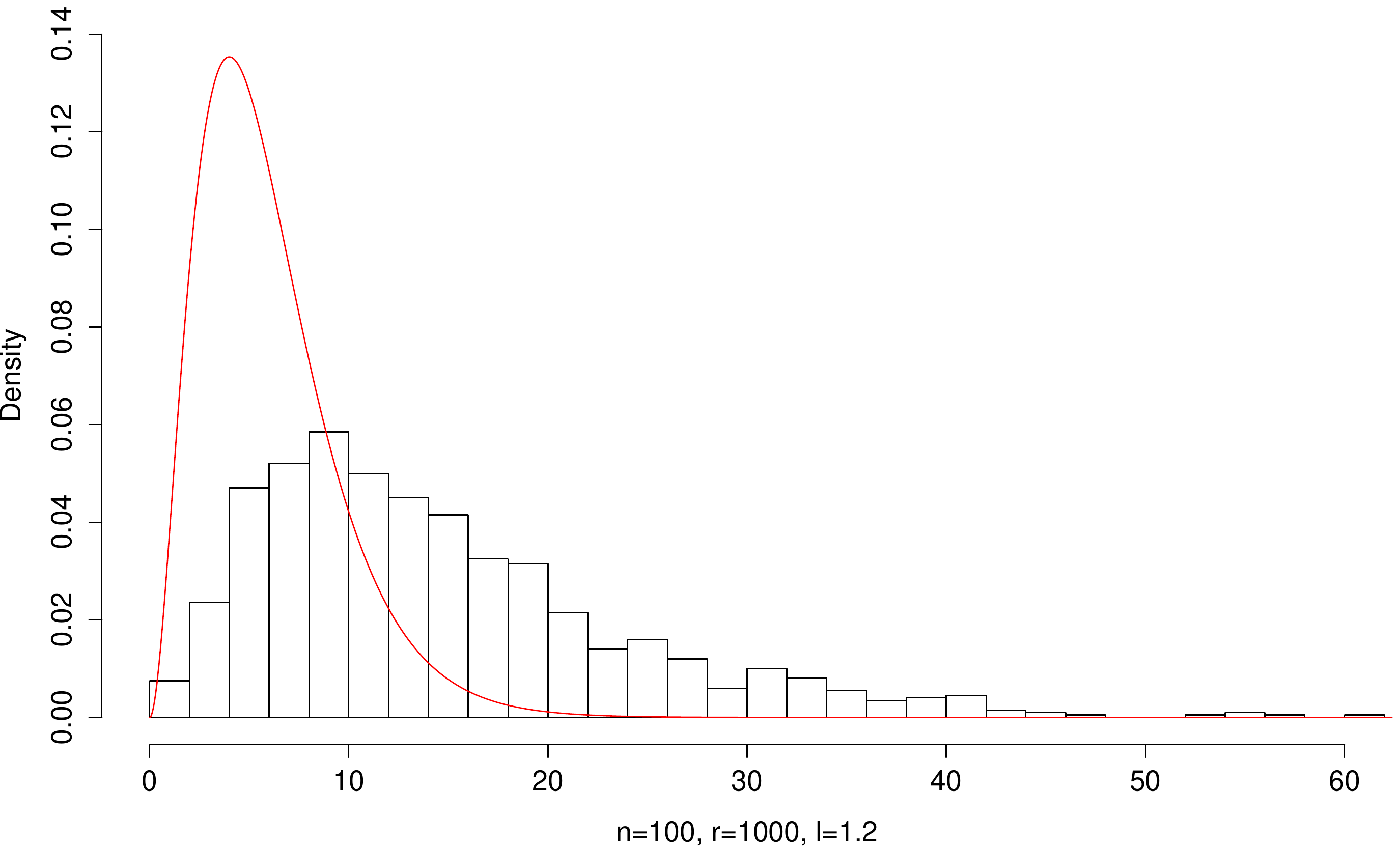}
      \caption{Histogram of Teststatistic with a tuned $\lambda_{msrL}=$l=1.2 for $n=100$}
      \label{fig2}
\end{figure}

First let us look at the question how the histogram of the teststatistic looks like for different $n$.
Here we use $p=500$ and $q=6$ with $J=(1,3,4,8,10,33)$, where the entries in $\beta_{J}$ are chosen randomly from a Uniform distribution on $[1,4]$. 
We also specify $\beta_{-J}$ to be the zero vector.
So the set $J$ that we are interested in, is in fact the same set as the active set of $\beta$. Furthermore $p-q$ gives the amout of sparsity.
Here we look at a sequence of $n=100,200,300,400,500,600,800$. As above, for each setting we calculate $1000$ simulations.
For each setting we plot the histogram for the teststatistic and compare it with the theoretical chi-squared distribution on $q=6$ degrees of freedom.
Figure \ref{fig1} and figure \ref{fig2} show the results.

The histograms show that with increasing $n$, we get a fast convergence to the true asymptotic chi-squared distribution. It is in fact true that we could
multiply the teststatistic with a constant $C_{n}\leq1$ in order to get the histogram match the chi-squared distribution. This reflects the theory.
Already with $n=400$ we get a very good approximation of the chi-quared distribution. But we see that the tuning of $\lambda_{msrL}$ is crucial for small $n$, see figure \ref{fig2}.

\begin{figure}[h!t]
  \centering
    \includegraphics[width=1\textwidth]{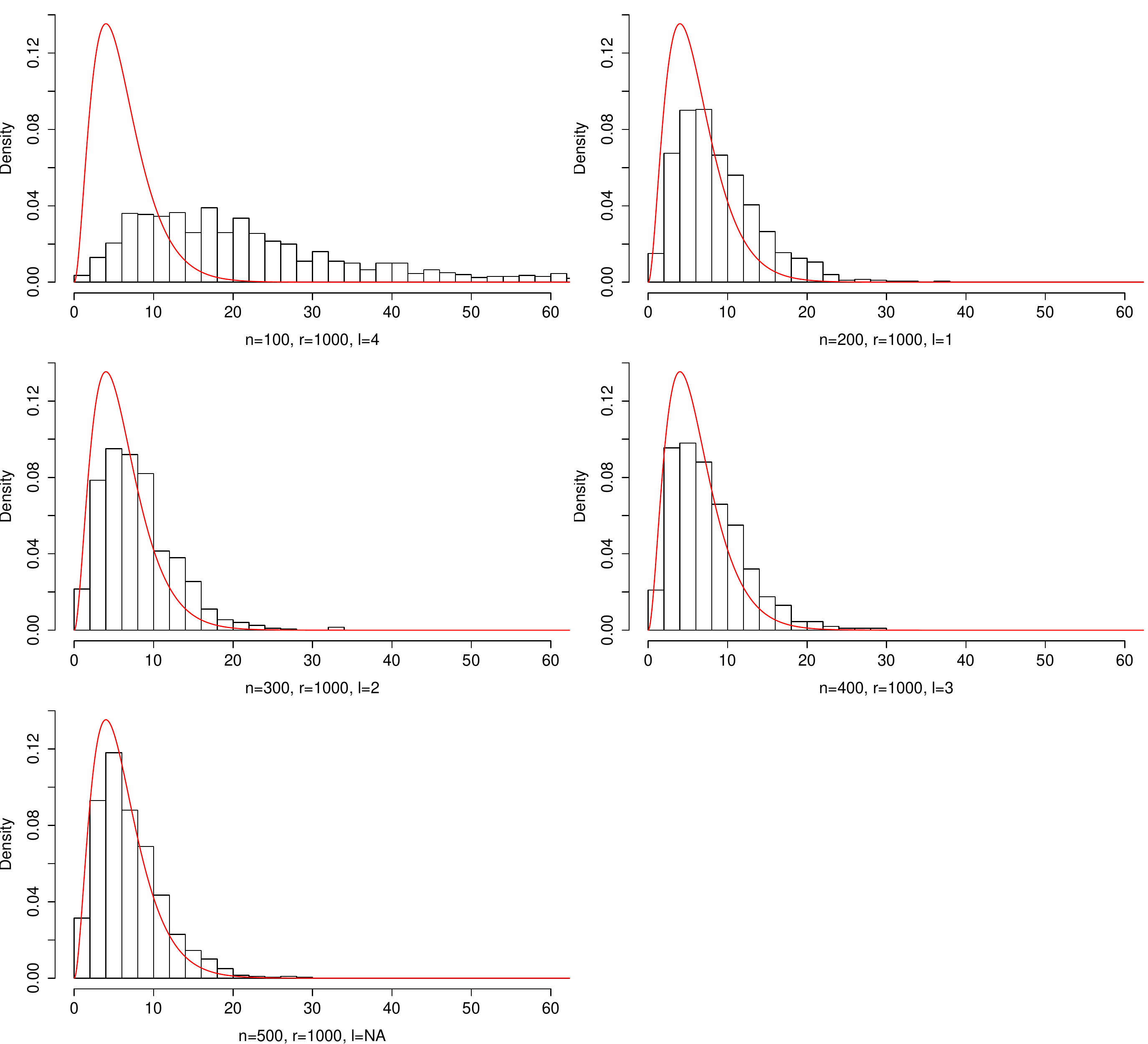}
      \caption{Histogram of Teststatistic, l=$\lambda$, the cross-validation lambda}
      \label{fig12}
\end{figure}

Next we try the same procedure but we are interested in what happens if we let $J$ and the active set $S_{0}$ not be the same set.
Here we take $J=(1,3,4,8,10,33)$ and $\beta_{0}$ is taken from the uniform distribution on $[1,4]$ on the set $S_{0}=(2,3,5,8,11,12,14,31)$. 
So only the indices $J\cap S_{0}= \{3,8\}$ coincide. Figure \ref{fig12} shows the results.

Compared to the case where $J$ and $S_{0}$ are the same set, it seems that this setting can handle small 
$n$ better than in the case where all the elements of $J$ are the nonzero indices of $\beta_{0}$.
So the previous case seems to be the harder case. Therefore we stick with $J=S_{0}=(1,3,4,8,10,33)$ for all the other simulations in Subsection \ref{sim:s1} and \ref{sim:s2}.

\subsection{Confidence level for an increasing $\lambda_{msrL}$}\label{sim:s1}
\begin{figure}[h!t]
    \centering
    \includegraphics[width=0.6\textwidth]{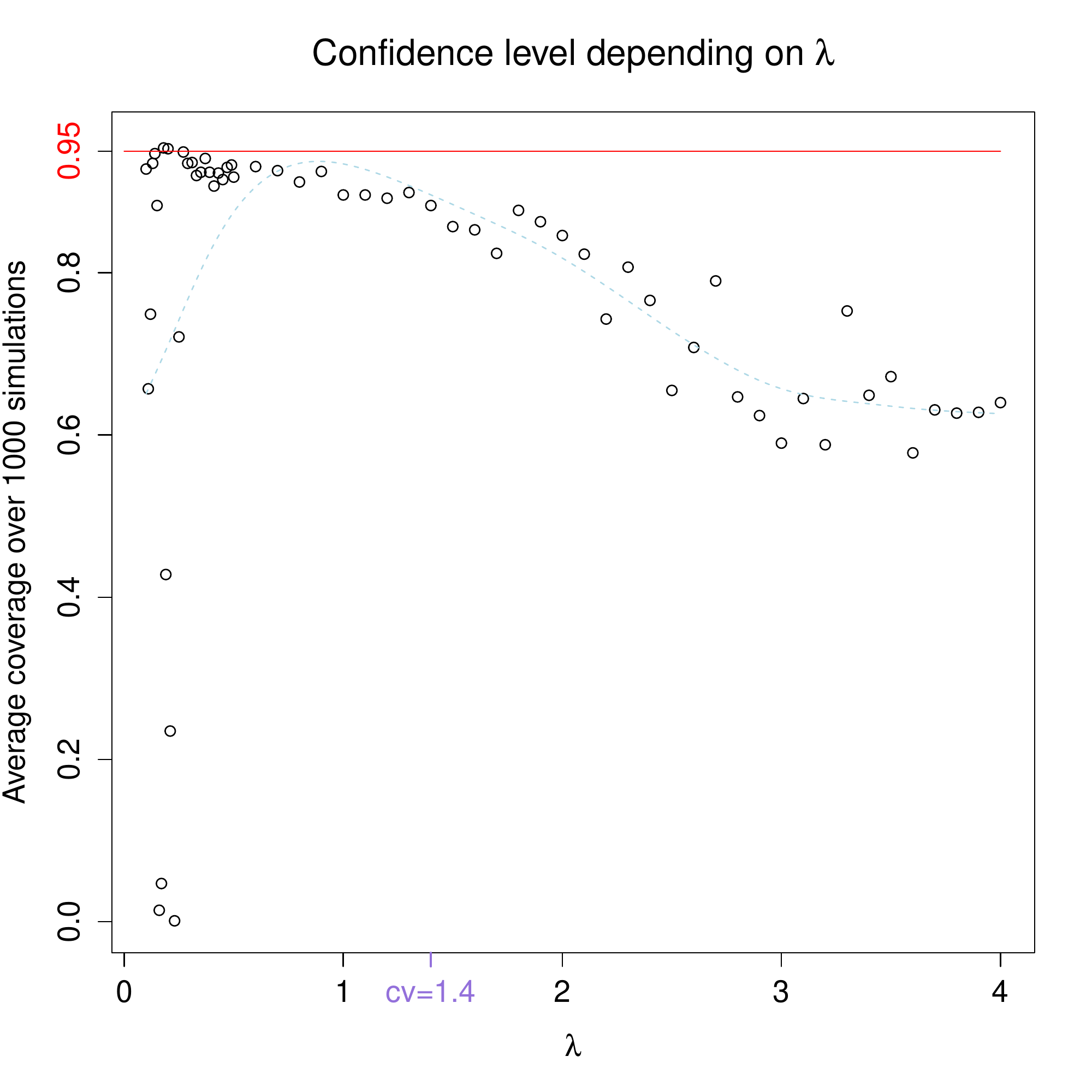}
    \caption{Average confidence level with fixed $n=400$ and $p=500$, increased $\lambda$}
    \label{fig5}
\end{figure}

Up until now we have not looked at the behaviour for different $\lambda$. We only used the cross-validation $\lambda$. 
So here we look at $n=400$, $p=500$ and we take $\lambda_{msrL}=(0.01,0.11,0.21,...,2.91)$ a fixed sequence.
Figure \ref{fig5} shows the results.

If we take $\lambda$ too low the behaviour breaks down. On the other hand, if $\lambda$ is too big,
we will not achieve a good average confidence level. The cross-validation $\lambda$ seems to be still a bit to high. So the cross-validation $\lambda$ could be better.

\subsection{Levelplot for n and p}\label{sim:s2}
\begin{figure}[h!t]
    \centering
    \includegraphics[width=1.1\textwidth, height=1\textwidth]{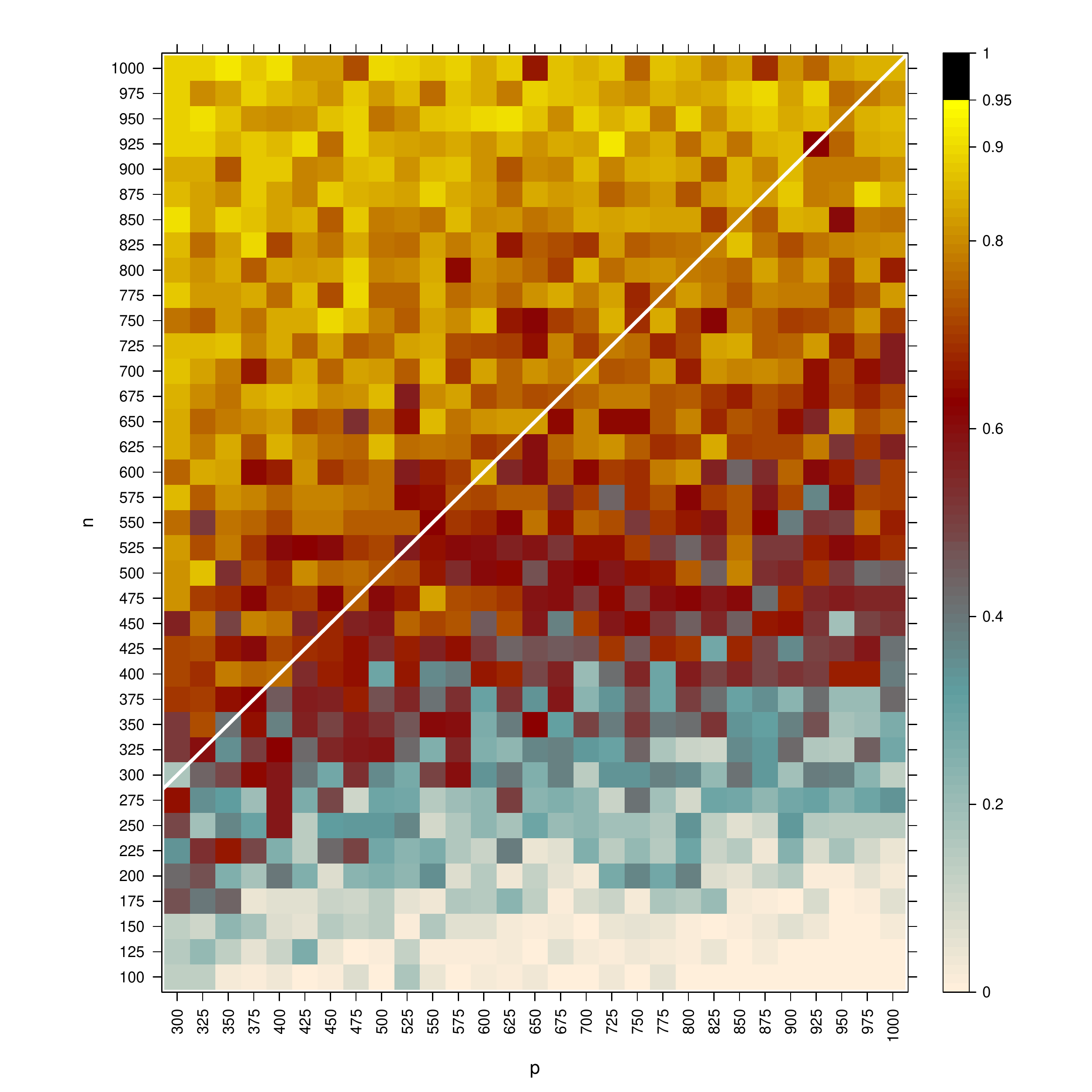}
    \caption{Levelplot of average confidence level for a grid $n,p$}
    \label{fig6}
\end{figure}

Not let us look at an overview of alot of different settings. We will use the levelplot to present the results
Here we use the cross-validation $\lambda$. We let $n$ and $p$ increase and look again at the average coverage of the confidence interval (average over the $1000$ simulations for each gridpoint).
The border between high and low dimensional cases is marked by the white line in figure \ref{fig6}.
Increasing $p$ does not worsen the procedure too much, which is very good. And, as expected, increasing the number of observations $n$ increases "the accuracy" of the average confidence interval.

\section{Discussion}\label{discussion.section}
We have presented a method for constructing confidence sets
for groups of variables which does not impose sparsity conditions on the input matrix
$X$. The idea is to use a loss function based on the nuclear norm of
the matrix of residuals. We called this the multivariate square-root Lasso as
it is an extension of the square-root Lasso in the multivariate case.

It is easy to see that when the groups are large, one needs the $\ell_2$-norm of the remainder
term $\| {\rm rem}\|^2 $ in Theorem \ref{main.theorem} to be of small order $\sqrt {|J|}$
in probability, using the representation
$\chi^2_{|J|}  = |J| + O_{\PP} (\sqrt {J} ) $. This leads to the requirement that
$ \sqrt n \lambda \| \hat \beta_{-J} - \beta_{j}^0 \|_1/\sigma_0  = o_{\PP} (1/ |J|^{1/4} )$,
i.e., that it
decreases faster for large groups. The paper \cite{Mitra} introduces a different
scheme for confidence sets, where there is no dependence on group size in the remainder term
after the normalization for large groups.
 Their idea is to use a group Lasso with a nuclear norm type of penalty on $\Gamma_J$ instead
 of the $\ell_1$-norm $\|\Gamma_J \|_1$ as we do in Theorem \ref{main.theorem}.
 Combining the approach of \cite{Mitra} with the result of Theorem \ref{structured.theorem} leads
 to a new remainder term which after normalization for large groups does not depend on group size {\it and} does not
 rely on sparsity assumptions on the design $X$. 
 
 The choice of the tuning parameter $\lambda$ for the construction used in Theorem \ref{main.theorem}
 is as yet an open problem. When one is willing to assume certain sparsity assumptions
 such that a bound for $\| \hat \beta - \beta^0\|_1$ is available, the tuning parameter can be
 chosen by trading off the size of the confidence set and the bias. When the rows
 of $X$ are i.i.d.\ random variables, a choice for $\lambda$ of order
 $\sqrt {\log p / n}$ is theoretically justified under certain conditions. Finally,
 smaller $\lambda$ give more conservative confidence intervals.
 Thus, increasing $\lambda$ will give one a ``solution path" of significant variables
 entering and exiting, where the number of ``significant" variables increases.
 If one aims at finding potentially important variables, one might want
 to choose a cut-off level here, i.e.\ choose $\lambda$ in such a way that
 the number of ``significant" variables is equal to a prescribed number.
  However, we have as yet no theory showing such a data-dependent
 choice of $\lambda$ is meaningful.
 
 A given value for $\lambda$ may yield sets which do not have the approximate
 coverage. These sets can nevertheless be viewed as giving a useful importance measure
 for the variables, an importance measure which avoids the possible  problems
 of other methods for accessing accuracy. For example,  when applied to all
 variables (after grouping) the confidence sets clearly also give results for the possibly weak variables.
 This is in contrast to post-model selection where the variables not selected are no
 longer under consideration.

\section{Proofs}\label{proofs.section} 

\subsection{Proof for the result for the multivariate square-root Lasso in Subsection \ref{multivariate.section}}

{\bf Proof of Lemma \ref{KKT.lemma}.} Let us write, for each $p \times q$ matrix $B$, the residuals as
  $ \Sigma (B):= (Y- X  B)^T (Y-X B)/n$.
  Let $\Sigma_{\rm min} (B)$ be the minimizer  of
  \begin{equation}\label{Sigma.equation}
  {\rm trace}  (\Sigma(B) \Sigma^{-1/2} ) + {\rm trace} ( \Sigma^{1/2}) 
  \end{equation}
 over $\Sigma$. Then  
 $ \Sigma_{\rm min} (B)$ equals $ \Sigma (B)$. To see this we invoke
  the reparametrization $\Omega:= \Sigma^{-1/2}$ so that $\Sigma^{1/2} = \Omega^{-1}$.
  We now minimize 
  $$ {\rm trace} (\Sigma(B) \Omega) + {\rm trace} ( \Omega^{-1}) $$
  over $\Omega >0$. 
  The matrix derivative with respect to $\Omega$ of 
  $ {\rm trace} (\Sigma(B) \Omega)$ is
  $ \Sigma (B)$.
   The matrix derivative of ${\rm trace} (\Omega^{-1}) $ with respect to
  $\Omega$ is equal to $-\Omega^{-2}$. Hence the minimizer $ \Omega_{\rm min} (B)$
  satisfies the equation
  $$\Sigma (B) - \Omega_{\rm min}^{-2} (B)=0 , $$
  giving 
  $$\Omega_{\rm min} (B)= \Sigma^{-1/2} (B).   $$
  so that
  $$ \Sigma_{\rm min} (B)= \Omega_{\rm min}^{-2} (B) = \Sigma (B) . $$
  Inserting this solution back in (\ref{Sigma.equation}) gives
  $2 {\rm trace} ( \Sigma^{1/2}  (B)) $ which is equal to $ 2 \| Y - X B \|_{\rm nuclear} /\sqrt n $.
  This proves the first part of the lemma.

  Let now for each $\Sigma >0$, $B (\Sigma)$ be the minimizer of 
  $$ {\rm trace}  (\Sigma(B) \Sigma^{-1/2} )+ 2\lambda_0 \| B \|_1 . $$
 By sub-differential calculus  we have 
 $$ X^T ( Y- X B) \Sigma^{-1/2} /n= \lambda_0 Z(\Sigma) $$
 where $\| Z ( \Sigma ) \|_{\infty} \le 1 $ and
 $Z_{k,j}(\Sigma )) = {\rm sign} ( B_{k,j} (\Sigma)) $ if $B_{k,j} (\Sigma ) \not= 0 $
 ($k=1 , \ldots , p$, $j=1 , \ldots q$). The KKT-conditions (\ref{multi-KKT.equation}) follow
 from $\hat B = B (\hat \Sigma)$. 
 \hfill $\sqcup \mkern -12mu \sqcap$

\subsection{Proof of the main result in Subsection \ref{main.section}}

{\bf Proof of Theorem \ref{main.theorem}.}
 We have
 $$M (\hat b_J -\hat \beta_J)
  =
 \hat T_J^{-1/2} (X_J - X_{-J} \hat \Gamma_J)^T   \epsilon /\sqrt n
 -  \hat T_J^{-1/2}(X_J - X_{-J} \hat \Gamma_J)^T   X( \hat \beta - \beta^0) /\sqrt n$$
 $$=  \hat T_J^{-1/2} (X_J - X_{-J} \hat \Gamma_J)^T   \epsilon /\sqrt n
 -  \hat T_J^{-1/2}(X_J - X_{-J} \hat \Gamma_J)^T   X_J( \hat \beta_J - \beta_J^0) /\sqrt n$$ $$
  -  \hat T_J^{-1/2}(X_J - X_{-J} \hat \Gamma_J)^T   X_{-J}( \hat \beta_{-J} - \beta_{-J}^0)/\sqrt n
 $$
 $$=    \hat T_J^{-1/2} (X_J - X_{-J} \hat \Gamma_J)^T   \epsilon /\sqrt n
 -M ( \hat \beta_J - \beta_J^0) 
 - \sqrt n \lambda \hat Z_J^T (\hat \beta_{-J}- \beta_{-J}^0 )$$
 where we invoked the KKT-conditions (\ref{gammaKKT.equation}).
 We thus arrive at 
  \begin{equation} \label{M2.equation}
  M (\hat b_J -\beta_J^0)=
   \hat T_J^{-1/2} (X_J - X_{-J} \hat \Gamma_J)^T   \epsilon /\sqrt n  +\sigma_0 {\rm rem},
    \end{equation}
    where
    $${\rm rem} = -\sqrt n \lambda \hat Z_J^T (\hat \beta_{-J}- \beta_{-J}^0 )/\sigma_0. $$
    The co-variance matrix of the first term 
     $\hat T_J^{-1/2} (X_J - X_{-J} \hat \Gamma_J)^T   \epsilon /\sqrt n $ in (\ref{M2.equation})
    is equal to
    $$\sigma_0^2 \hat T_J^{-1/2} (X_J - X_{-J} \hat \Gamma_J)^T  (X_J - X_{-J} \hat \Gamma_J)\hat T_J^{-1/2}/n=
    \sigma_0^2 I $$
    where $I$ is the identity matrix with dimensions $|J| \times |J|$. It follows that
    this term is $|J|$-dimensional standard normal scaled with $\sigma_0$.
The remainder term can be bounded using the dual norm inequality for each
entry:
$$|{\rm rem}_j| \le \sqrt n \lambda  \max_{k \notin J} |(\hat Z_J)_{k,j} | \| \hat \beta_{-J} - \beta_{-J}^0 \|_1/\sigma_0 \le
\sqrt n \lambda \| \hat \beta_{-J} - \beta_{-J}^0 \|_1/\sigma_0  $$
since by the KKT-conditions (\ref{gammaKKT.equation}), we have $\| \hat Z_J \|_{\infty} \le 1 $.
\hfill $\sqcup \mkern -12mu \sqcap$

\subsection{Proofs of the theoretical result for the square-root Lasso in Section \ref{theory.section}}

{\bf Proof of Lemma \ref{Gauss.lemma}.} Without loss of generality we can assume $\sigma_0^2 =1$.
From \cite{laurent2000adaptive} we know that for all $t>0$
$$ \PP \biggl( \| \epsilon \|_n^2 \le 1- 2 \sqrt {t/n} \biggl) \le \exp [-t]  $$
and
$$\PP\biggl ( \| \epsilon \|_n^2 \ge 1+ 2 \sqrt {t/n} + 2t/n \biggl) \le \exp [-t] .$$

Apply this with $t= \log(1/{\underline \alpha})$ and $t = \log (1/\bar \alpha)$
respectively.
Moreover $X_j^T \epsilon / n \sim {\cal N} (0, 1/n )$ for all $j$.
Hence for all $t>0$
$$ \PP \biggl( |X_j^T \epsilon| / n\ge \sqrt {2t / n} \biggl) \le 2 \exp[-t] , \ \forall \ j . $$
It follows that
$$\PP \biggl ( \| X^T \epsilon \|_{\infty}/ n \ge \sqrt {2(t + \log (2p) ) / n}\biggl ) \le  \exp[-t] . $$
\hfill $\sqcup \mkern -12mu \sqcap$ 

 {\bf Proof of Lemma \ref{sigma.lemma}.} Suppose $\hat R \le R$ and $\| \epsilon \|_n \ge {\underline \sigma}$.
    First we note that the inequality
    (\ref{eta.equation}) gives
    $$\lambda_0 \| \beta^0 \|_1 / \| \epsilon \|_n \le 2 \biggl ( \sqrt { 1 + (\eta/2)^2} -1 \biggr ) . $$
    For the upper bound for $\| \hat \epsilon \|_n$ we use that
    $$\| \hat \epsilon \|_n + \lambda_0 \| \hat \beta \|_1 \le \| \epsilon \|_n + \lambda_0 \| \beta^0 \|_1 $$
    by the definition of the estimator. 
    Hence
    $$\| \hat \epsilon \|_n \le \| \epsilon \|_n + \lambda_0 \|\beta^0 \|_1 \le \biggl [ 1+ 
    2 \biggl ( \sqrt { 1 + (\eta/2)^2} -1 \biggr ) \biggr ] \| \epsilon\|_n \le (1+ \eta) \| \epsilon \|_n . $$
    For the lower bound for $\| \hat \epsilon \|_n$ we use the convexity of both the loss
    function and the penalty. 
     Define
     $$\alpha := { \eta \| \epsilon \|_n \over \eta \| \epsilon\|_n + \| X (\hat \beta - \beta^0 ) \|_n} .$$
     Note that $0 < \alpha \le 1$. 
     Let $\hat \beta_{\alpha}$ be the convex combination $\hat \beta_{\alpha} :=  \alpha  \hat \beta + (1-\alpha) \beta^0 $. 
     Then
     $$ \| X ( \hat \beta_{\alpha} - \beta^0) \|_n = \alpha \| X ( \hat \beta - \beta^0) \|_n =
     {\eta \| \epsilon \|_n  \| X ( \hat \beta - \beta^0) \|_n \over \eta \| \epsilon \|_n  + \| X ( \hat \beta - \beta^0) \|_n  }
     \le \eta \| \epsilon \|_n . $$
     Define $\hat \epsilon_{\alpha} := Y - X \hat \beta_{\alpha}  $. 
     Then, by convexity of $\| \cdot \|_n$ and $\| \cdot \|_1 $,
     $$\| \hat \epsilon_{\alpha} \|_n + \lambda_0 \| \hat \beta_{\alpha} \|_1 \le \alpha  \| \hat \epsilon \|_n +
     \alpha \lambda_0 \| \hat \beta \|_1 + (1- \alpha) \| \epsilon \|_n + (1- \alpha ) \lambda_0 \| \beta^0 \|_1 $$
     $$ \le \| \epsilon \|_n + \lambda_0 \| \beta^0 \|_1 $$
     where in the last step we again used that $ \hat \beta$ minimizes $\| Y - X \beta \|_n + 
     \lambda_0 \| \beta \|_1$. 
     Taking squares on both sides gives
     \begin{equation} \label{squares.equation}
     \| \hat \epsilon_{\alpha} \|_n^2 + 2 \lambda_0 \| \hat \beta_{\alpha} \|_1  \| \hat \epsilon_{\alpha } \|_n+
     \lambda_0^2 \| \hat \beta_{\alpha}  \|_1^2 \le 
     \| \epsilon \|_n^2 + 2 \lambda_0  \| \beta^0 \|_1 \| \epsilon \|_n +
     \lambda_0^2 \| \beta^0 \|_1^2  . 
     \end{equation}
     But
     $$ \| \hat \epsilon_{{\alpha }} \|_n^2 = \| \epsilon \|_n^2 - 2 \epsilon^T X( \hat \beta_{\alpha} - 
     \beta^0 ) /n + \| X ( \hat \beta_{\alpha} -\beta^0 ) \|_n^2 $$
     $$ \ge \| \epsilon \|_n^2 - 2 R \| \hat \beta_{\alpha} - \beta^0 \|_1 \| \epsilon \|_n + \| X ( \hat \beta_{\alpha} -\beta^0 ) \|_n^2 
      $$
      $$ \ge \| \epsilon \|_n^2 - 2 R \| \hat \beta_{\alpha} \|_1\| \epsilon \|_n  -2R \|  \beta^0 \|_1 \| \epsilon \|_n + \| X ( \hat \beta_{\alpha} -\beta^0 ) \|_n^2 . $$
     Moreover, by the triangle inequality
     $$ \| \hat \epsilon_{\alpha} \|_n \ge \| \epsilon \|_n - \| X ( \hat \beta_{\alpha} - \beta^0) \|_n 
     \ge (1- \eta ) \| \epsilon \|_n . $$
     Inserting these two inequalities into (\ref{squares.equation}) gives
     $$ \| \epsilon \|_n^2 - 2 R \| \hat \beta_{\alpha}\|_1 \| \epsilon \|_1  - 2 R \| \beta^0 \|_1 \| \epsilon \|_n + \| X ( \hat \beta_{\alpha} -\beta^0 ) \|_n^2      +  2 \lambda_0 (1- \eta) \| \hat \beta_{\alpha} \|_1 \| \epsilon \|_n  +
     \lambda_0^2 \| \hat \beta_{\alpha}  \|_1^2 $$ $$\le   \| \epsilon \|_n^2 + 2 \lambda_0  \| \beta^0 \|_1 \| \epsilon \|_n+
     \lambda_0^2 \| \beta^0 \|_1^2  $$
     which implies by the assumption $\lambda_0 (1- \eta ) \ge R$
    $$  \| X ( \hat \beta_{\alpha} -\beta^0 ) \|_n^2\le 2 (\lambda_0 + R) \| \beta^0 \|_1 \| \epsilon \|_1 +
    \lambda_0^2 \| \beta^0 \|_1^2 . $$
    $$ \le 4 \lambda_0 \| \beta^0 \|_1 \| \epsilon \|_1 +  \lambda_0^2 \| \beta^0 \|_1^2 $$
    where in the last inequality we used $R \le (1- \eta ) \lambda_0 \le \lambda_0 $.
    But continuing we see that we can write the last expression as
    $$  4 \lambda_0 \| \beta^0 \|_1 \| \epsilon \|_1 +  \lambda_0^2 \| \beta^0 \|_1^2 =
   \biggl ( ( \lambda_0 \| \beta_0 \|_1/ \| \epsilon_n \|_n + 2)^2 - 4 \biggr ) \| \epsilon \|_n^2.
      $$
      Again invoke the $\ell_1$-sparsity condition
      $$ \lambda_0 \| \beta^0 \|_1 / \| \epsilon \|_n \le 2 \biggl ( \sqrt { 1 + (\eta/2)^2} -1 \biggr )  $$
    to get
   $$ \biggl ( ( \lambda_0 \| \beta_0 \|_1/ \| \epsilon_n \|_n + 2)^2 - 4 \biggr ) \| \epsilon \|_n^2 \le
   { \eta^2  \over 4 } \| \epsilon \|_n^2 . $$
      We thus established that     
      $$ \| X ( \hat \beta_{\alpha} -\beta^0 ) \|_n \le { \eta \| \epsilon \|_n \over 2} . $$
     Rewrite this to
     $$ {\eta \| \epsilon \|_n  \| X ( \hat \beta - \beta^0) \|_n \over \eta \| \epsilon \|_n  + \| X ( \hat \beta - \beta^0) \|_n  }
     \le { \eta \| \epsilon \|_n \over 2}, $$
     and rewrite this in turn to
     $$\eta \| \epsilon \|_n  \| X ( \hat \beta - \beta^0) \|_n \le
    {  \eta^2 \| \epsilon \|_n^2 \over 2 } + { \eta \| \epsilon \|_n \| X ( \hat \beta - \beta^0 ) \|_n \over 2} $$
    or 
    $$\| X ( \hat \beta - \beta^0 ) \|_n \le \eta \| \epsilon \|_n .$$
    But then, by repeating the argument, also
    $$\| \hat \epsilon \|_n \ge \| \epsilon \|_n - \| X ( \hat \beta - \beta^0 ) \|_n \ge (1- \eta ) \| \epsilon \|_n . $$
    \hfill $\sqcup \mkern -12mu \sqcap$

{\bf Proof of Theorem \ref{oracle.theorem}.}  
Throughout the proof we suppose
$\hat R \le R$ and $\| \epsilon \|_n \ge {\underline \sigma}$.
Define the Gram matrix $\hat \Sigma := X^T X/ n $. 
Let $\beta \in \R^p$ and $S:=S_{\beta}= \{ j :\ \beta_j \not= 0 \}$. 
    If 
  $$ ( \hat \beta - \beta)^T \hat \Sigma  ( \hat \beta - \beta^0 )  \le -\delta \underline{\lambda} \|  \hat \beta - \beta \|_1 \| \epsilon \|_n $$  we find 
    $$2\delta \underline{\lambda}\|  \hat \beta - \beta \|_1 \|  \epsilon \|_n  +\| X ( \hat \beta - \beta^0 )\|_n^2  
  $$ $$ =2\delta \underline{\lambda}\| \hat \beta - \beta \|_1  \| \epsilon \|_n + \| X (\beta - \beta^0 ) \|_n^2 
      - \| X ( \beta - \hat \beta ) \|_n^2  + 2 ( \hat \beta - \beta)^T \hat \Sigma  ( \hat \beta - \beta^0 ) $$ $$
    \le  \| X (\beta - \beta^0 ) \|_n^2 .$$
   So then
    we are done.
  
  Suppose now that
  $$(\hat \beta - \beta)^T \hat \Sigma ( \hat \beta -\beta^0 ) \ge -\delta \underline{\lambda} \|  \hat \beta - \beta \|_1 \| \epsilon \|_n 
  .$$
  By the KKT-conditions (\ref{KKT.equation})
  $$(\hat \beta - \beta )^T \hat \Sigma (\hat \beta -\beta^0 )  + \lambda_0   \| \hat \beta \|_1 
  \| \hat \epsilon \|_n  \le \epsilon^T X(\hat \beta - \beta)/n + \lambda_0 \|\beta \|_1 \| \hat \epsilon \|_n  . $$
  By the dual norm inequality and since $\hat R \le R$
  $$| \epsilon^T X(\hat \beta - \beta)| /n   \le R \| \hat \beta - \beta \|_1 \| \epsilon \|_n 
    . $$
  Thus
  $$ (\hat \beta - \beta )^T\hat \Sigma (\hat \beta -\beta^0 )  + \lambda_0  \| \hat \beta\|_1\| \hat \epsilon \|_n   \le 
  R \| \hat \beta - \beta \|_1 \| \epsilon \|_n + \lambda_0 \| \beta\|_1 \| \hat \epsilon \|_n . $$
 This implies by the triangle inequality
 $$
  (\hat \beta - \beta )^T\hat \Sigma (\hat \beta -\beta^0 )  + (\lambda_0 \| \hat \epsilon \|_n - R \| \epsilon \|_n )  \| \hat \beta_{-S} \|_1 \le 
   (\lambda_0 \| \hat \epsilon \|_n +  R \| \epsilon \|_n  ) \| \hat \beta_S - \beta \|_1  . 
   $$ 
 We invoke the result of  Lemma \ref{sigma.lemma} which says that that $(1- \eta ) \| \epsilon \|_n \le \| \hat \epsilon \|_n \le (1+ \eta ) \| \epsilon \|_n $. This gives
   \begin{equation}\label{tea1.equation}
  (\hat \beta - \beta )^T\hat \Sigma (\hat \beta -\beta^0 )  + \underline{\lambda}   \| \hat \beta_{-S} \|_1 \| \epsilon \|_n\le 
   (\lambda_0 (1+ \eta) +  R)    \| \hat \beta_S - \beta \|_1   \| \epsilon \|_n. 
   \end{equation} 
 Since $  (\hat \beta - \beta )^T \hat \Sigma (\hat \beta -\beta^0 ) \ge - \delta \underline{\lambda}
   \|  \epsilon \|_n  \| \hat \beta - \beta \|_1$ this gives
 $$(1-\delta) \underline{\lambda} \|  \hat \beta_{-S} \|_1 \| \epsilon \|_n \le (\lambda_0 (1+ \eta) + R+
 \delta \underline{\lambda}) \| \hat \beta_S - \beta \|_1 \| \epsilon \|_n  
 = \bar \lambda \| \hat \beta_S - \beta \|_1 \| \epsilon \|_n . $$
 or
 $$ \|  \hat \beta_{-S} \|_1 \le  L \| \hat \beta_S - \beta \|_1.  $$
 But then
 \begin{equation}\label{coffee1.equation}
 \| \hat \beta_S - \beta \|_1  \le \sqrt {|S|}   \| X ( \hat \beta - \beta ) \|_n / \hat \phi(L,S) .
 \end{equation}
 Continue with inequality (\ref{tea1.equation}) and apply the inequality $ab \le (a^2 + b^2 )/2$
 which holds for all real valued $a$ and $b$:
 $$(\hat \beta - \beta )\hat \Sigma  (\hat \beta -\beta^0 )  + \underline{\lambda}  \| \hat \beta_{-S}\|_1 \| \epsilon \|_n+ \delta \underline{\lambda}
 \| \hat \beta_S - \beta \|_1 \| \epsilon \|_n $$ $$ \le 
   \bar \lambda \| \epsilon\|_n \sqrt {|S|}  \| X ( \hat \beta - \beta ) \|_n  /
   \hat \phi (L,S) $$
   $$ \le {1 \over 2} {\bar \lambda}^2  { |S| \| \epsilon \|_n^2 \over  \hat \phi^2 (L,S) }  + {1 \over 2} \| X ( \hat \beta - \beta ) \|_n^2 . $$
   Since    $$
   2 (\hat \beta - \beta )^T \hat \Sigma ( \hat \beta - \beta^0 )= \| X ( \hat \beta - \beta^0 ) \|_n^2  - \| X (\beta - \beta^0 ) \|_n^2 +
   \| X ( \hat \beta - \beta ) \|_n^2  ,$$
   we obtain
   $$\| X (\hat \beta - \beta^0 ) \|_n^2 + 2\underline{\lambda}  \| \hat \beta_{-S}\|_1\| \epsilon \|_n   +  2 \delta \underline{\lambda}  \|  \hat \beta_S - \beta\|_1 \| \epsilon \|_n  $$ $$
\le \| X (\beta - \beta^0 ) \|_n^2 + {\bar \lambda}^2  |S|\| \epsilon \|_n^2  /\hat \phi^2 (L,S)  . $$
  \hfill $\sqcup \mkern -12mu \sqcap$

\subsection{Proofs of the illustration assuming (weak) sparsity in Section \ref{lr.section}}

{\bf Proof of Lemma \ref{lr.lemma}.}
Define $\lambda_*:= {\bar \lambda} \| \epsilon \|_n / \hat \Lambda_{\rm max}  (S_0) $
and for $j=1 , \ldots , p $, 
$$\beta_j^* = \beta_j^0 {\rm l} \{ | \beta_j^0 | > \lambda_*  \} . $$
  Then   
  $$\| X ( \beta^* - \beta^0 ) \|_n^2 \le \hat \Lambda_{\rm max}^2 (S_0) \| \beta^* - \beta^0 \|_2^2  
  \le \hat \Lambda_{\rm max}^2 (S_0) \lambda_*^{2- r} \rho_r^r  $$ $$=
  \bar \lambda^{2-r}  \hat \Lambda_{\rm max}^r (S_0) \rho_r^r \| \epsilon \|_n^{2-r} \le
 \bar \lambda^{2-r}  \hat \Lambda_{\rm max}^r (S_0) \rho_r^r \| \epsilon \|_n^{2-r} / \hat \phi^2 (L, S_*) $$
 where in the last inequality we used $\hat\phi (L, S_*) \le 1$. 
  Moreover, noting that
  $S_{\beta^*} =\hat S_*  =\{ j :\ | \beta_j^0 | > \lambda^* \} $ we get
  $$| S_{\beta^*} | \le \lambda_*^{-r} \rho_r^r = \bar \lambda^{-r} \| \epsilon \|_n^{-r} \hat \Lambda_{\rm max}^r (S_0) 
  .$$
 Thus
  $$\bar \lambda^2 | S_{\beta^*} | \| \epsilon \|_n^2  /\hat \phi^{2} (L, S_{\beta^*} ) \le
  \bar \lambda^{2-r}  \hat \Lambda_{\rm max}^r (S_0)\rho_r^r  \| \epsilon \|_n^{2-r}/ \hat \phi^{2} (L, \hat S_* ) .
  $$
  Moreover
  $$\| \beta^* - \beta_0 \|_1 \le \lambda_*^{1-r } \rho_r^r = {\bar \lambda}^{1-r}  \| \epsilon \|_n^{1-r}   \hat  \Lambda_{\rm max}^r  (S_0) / \hat \phi^{2}  (L, \hat S_*)  ,$$
  since $\hat \phi^2 (L, \hat S_*) / \hat \Lambda_{\rm max} (S_0) \le 1 $. 
  \hfill $\sqcup \mkern -12mu \sqcap$
  
  {\bf Proof of Lemma \ref{bound.lemma}.} The $\ell_1$-sparsity condition (\ref{eta.equation}) holds with $\eta \le 1/3$.
 Theorem \ref{oracle.theorem} with $\lambda_0(1- \eta)  = 2R$ gives
 ${\underline \lambda} = \lambda_0 (1- \eta) -R = R$ and
 $3 R \le \bar \lambda = \lambda_0 (1+ \eta) +R + \delta {\underline \lambda} \le
 (5+ \delta ) R $. We take $\delta = 1/7$.
 Then $L= \bar \lambda/ ((1- \delta) {\underline \lambda})\le (5+ \delta )/ (1-\delta) = 6$. 
 Set $\hat S_* := \{ j : \ | \beta_j^0 | > \bar \lambda \| \epsilon \|_n / \hat \Lambda_{\rm max} (S_0)\} $.
 On the set where $\| \epsilon \|_n \ge {\underline \sigma}$ we have
 $\hat S_* \subset S_*$ since $\bar \lambda \ge 3 R$.
 We also have $\bar \lambda / (\delta {\underline \lambda} ) \le 6^2 $. 
 Hence, using the arguments of Lemma \ref{lr.lemma} and the
 result of Theorem \ref{oracle.theorem}, we get on the set $\hat R \le R$ and
 $\| \epsilon \|_n \ge {\underline \sigma}$, 
 $$   { \| \hat \beta - \beta^0 \|_1 \over \| \epsilon \|_n }  \le 
 {\bar \lambda}^{1-r} \biggl (  1   +
 {6^2  \hat \Lambda_{\rm max}^r (S_0) \over  \hat \phi^{2}  (6, S_*)}
\biggr ) \biggr ( { \rho_r  \over \| \epsilon \|_n   } \biggr )^r  .$$
Again, we can bound here $1/ \| \epsilon \|_n$ by $1/ {\underline \sigma} $.  
We can moreover bound $\bar \lambda $ by $6R$.  
Next we see that on the set where $\hat R \le R$ and
 $\| \epsilon \|_n \ge {\underline \sigma}$, by Lemma \ref{sigma.lemma},
$$ \hat \sigma \ge (1- \eta) \| \epsilon \|_n \ge (1- \eta) {\underline \sigma} . $$
The $\ell_0$-bound follows in the same way, inserting $\beta = \beta^0$ in Theorem \ref{main.theorem}.
Invoke Lemma \ref{Gauss.lemma} to show that the set $\{ \hat R \le R \cap \| \epsilon \|_n \ge 
\underline \sigma \}$ has probability at least $1- \alpha_0 - {\underline \alpha} $. 
\hfill $\sqcup \mkern -12mu \sqcap$

\subsection{Proof of the extension to structured sparsity in Section \ref{structured.section}} 

{\bf Proof of Theorem \ref{structured.theorem}.} This follows from exactly the same arguments as used in the proof
of Theorem \ref{main.theorem} as
the KKT-conditions (\ref{gammaKKT.equation}) with general norm $\Omega$ imply that
 $$
 \| X_{-J}^T ( X_J - X_{-J} \hat \Gamma_J ) \hat T_J^{-1/2}/n \|_{\infty, \Omega_* } \le  \lambda  .
$$
\hfill $\sqcup \mkern -12mu \sqcap$

\bibliographystyle{plainnat}
\bibliography{reference}


\end{document}